\documentclass[12pt]{article}
\usepackage[centertags]{amsmath}
\usepackage{amsfonts}
\usepackage{amssymb}
\usepackage{amsthm}
\usepackage{newlfont}
\usepackage[dvips]{graphicx}
 \newtheorem{thm}{Theorem}[section]
 
 \newtheorem{lem}[thm]{Lemma}
 
 \newtheorem{exmp}[thm]{Example}
 
 \theoremstyle{definition}
 \newtheorem{defn}{Definition}[section]

 \newtheorem{rem}[thm]{Remark}
 \numberwithin{equation}{section}

\title{On the practical global  uniform asymptotic stability of stochastic differential equations\thanks{The research of T. Caraballo has been partially supported by FEDER and the Ministerio de Econom\'{i}a y Competitividad (Spain) under grant MTM2011-22411, and Junta de Andaluc\'{i}a (Spain) under Proyecto de Excelencia P12-FQM-1492 and the Ayudas de consolidaci\'{o}n for the research group FQM314.}}

\vspace{8pt}

\author{Tom\'as Caraballo$^b$   \ \ \  Mohamed Ali Hammami$^a$  \ \ \ Lassaad Mchiri$^a$  \\
\\
{\small  $^a$ University of Sfax, Faculty of Sciences of Sfax, Department of Mathematics,  Tunisia}\\
{\small   E.mail: MohamedAli.Hammami@fss.rnu.tn}\\
{}
\\
 {\small  $^b$  Universidad de Sevilla,}\\
  {\small   Depto. Ecuaciones Diferenciales y An\'alisis Num\'erico,}\\
 {\small   Facultad de Matem\'aticas, }\\
  {\small  Apdo. de Correos 1160}\\
  {\small  41080-Sevilla (Spain)}\\
 }
\date {}
\begin{document}\maketitle

\vspace{10pt}
\begin{abstract}
The method of Lyapunov functions is one of the most effective ones for the investigation of
stability of dynamical systems, in particular, of stochastic differential systems. The main purpose
of the paper is the analysis of the stability of stochastic differential equations by using Lyapunov
functions  when the origin is not necessarily an equilibrium point. The global uniform boundedness and the global practical uniform exponential stability of solutions of  stochastic differential equations based on  Lyapunov techniques are investigated. Furthermore, an example is given to illustrate the applicability of the main result.
\end{abstract}
\vspace{10pt}


\section{Introduction}

The main aim of this paper is to establish some criteria for the global practical uniform exponential stability of a class of nonlinear stochastic differential equations of the form 
\begin{equation}\label{e1}
dx(t)=f(x(t),t)dt+g(x(t),t)dW(t).
\end{equation}
 As it is well known, stochastic differential equations are useful for modeling physical, technical, biological and economical dynamical systems in which significant uncertainly is present.
 
To investigate stability for stochastic differential equations, we usually study the stability and attractivity of solutions in a
neighborhood of the origin considered as an equilibrium point.  The equations are worth being investigated since the solutions with different initial values have similar large-time properties. When $f(0,t)=g(0,t)=0$, the stability and the exponential stability of equation (\ref{e1}) have received very much attention in the literature. We would like to mention here the references \cite{5}-\cite{12}, among others.

Several interesting and important variants to Lyapunov's original concepts of practical stability were proposed in \cite{1}-\cite{4}. When the origin is not necessarily an equilibrium point,
we can study the asymptotic stability of solutions with respect to a small neighborhood of the origin.
The goal is to analyze systems whose desired behavior is asymptotic stability about a small ball of the state space or a close approximation to this. Almost all the state trajectories are bounded and approach a sufficiently small neighborhood of the origin. One also wishes that the state approaches the origin (or some sufficiently small neighborhood of it) in a sufficiently fast manner, for instance, exponentially. This property is referred to as practical stability.

Therefore, in this paper we state some sufficient conditions to ensure the almost sure practical global uniform asymptotic stability and the exponential stability of the solution of equation (\ref{e1}).
The content of the paper is as follows. In Section $2$, we introduce some auxiliary facts and results. In Section $3$, we prove the main theorem about the global practical uniform exponential stability of stochastic differential equation and exhibit an example as an application of this study.

\section{Auxiliary Facts and Results}

 Consider the following $d$-dimensional stochastic differential equation (SDE):
\begin{equation}\label{la}
dx(t)=f(x(t),t)dt+g(x(t),t)dW_{t}, \quad \forall t\geq0,
\end{equation}
where $f:\mathbb{R}^{d}\times\mathbb{R}_{+}\longrightarrow\mathbb{R}^{d}$ and $g:\mathbb{R}^{d}\times\mathbb{R}_{+}\longrightarrow\mathbb{R}^{d\times m}$, $x(t)=(x_{1}(t),...,x_{d}(t))^{T}\in\mathbb{R}^{d}$ and $W_{t}=(W_{1}(t),....,W_{m}(t))^{T}$ is an $m$-dimensional Brownian motion defined on a complete probability space $(\Omega,\mathcal{F},\mathbb{P})$.
The Euclidian norm of a vector $x\in\mathbb{R}^{d}$ is denoted by $|x|$.

We assume that both $f$ and $g$ satisfy the following conditions:
\begin{equation}
|f(x,t)|\leq C\big(1+|x|\big); \quad |g(x,t)|\leq C\big(1+|x|\big),\ \textrm{for all}\ t\geq0, x\in\mathbb{R}^{d},
\end{equation}
\begin{equation}
|f(x,t)-f(y,t)|\leq C|x-y|; \quad|g(x,t)-g(y,t)|\leq C|x-y|, \ \textrm{for all}\ t\geq0, x,y\in\mathbb{R}^{d}.
\end{equation}
Hence, for any given initial value $x_{0}\in\mathbb{R}^{d}$, the SDE has a unique global solution denoted by $x(t,t_{0},x_{0})$, such that $x(t_0)=x_0$, and defined for all $t\geq t_{0}\geq0$. In what follows we use $x(t,t_{0},x_{0})$, or simply $x(t)$ if no confusion is possible, to denote the solution of (\ref{la})  with initial value $x_0$ at the initial time $t_0$.\\
We suppose that  there exists $t$ such that $f(0,t)\neq0$ and $g(0,t)\neq0$. In this case, the SDE does not possess the trivial  solution $x(t,t_0,0)=0$. Therefore, we will study the asymptotic stability of the SDE when $0$ in not an equilibrium point, but in a small neighborhood of the origin in terms of convergence of solution in probability to a small ball $B_{r}:=\{x\in\mathbb{R}^{d}:|x|\leq r\}$, $r>0$.

\begin{defn}
 The solution of (\ref{la}) is said to be globally uniformly bounded in probability if, for each $\alpha>0$, there exists $ c=c(\alpha)>0$ (independent of $t_{0}$) such that,
\begin{equation}
\textrm{for every} \ t_{0}\geq0,\ \textrm{and all} \ x_0\in\mathbb{R}^d \ \textrm{with}\ \quad |x_{0}|\leq\alpha,\quad \sup_{t\geq t_{0}}|x(t,t_{0},x_{0})|\leq c(\alpha),  \quad \textrm{a.s.}
\end{equation}
\end{defn}
\begin{defn}
$(i)$ The ball $B_{r}$ is said to be uniformly stable in probability if,
for each $\epsilon\in]0,1[$ and $k>r$, there exists $\delta=\delta(\epsilon,k)>0$ such that
\begin{equation}
\mathbb{P}\big(|x(t,t_{0},x_{0})|<k, \forall t\geq t_{0}\geq0\big)\geq1-\epsilon\quad \text{for all}\quad |x_{0}|<\delta.
\end{equation}
$(ii)$ The ball $B_{r}$ is said to be globally uniformly stable in probability if it is uniformly stable in probability and the solution of (\ref{la}) is globally uniformly bounded in probability.
\end{defn}
\begin{defn}
The ball $B_{r}$ is said to be globally uniformly  attractive in probability if for each
$\epsilon\in]0,1[$, $k>r$ and $c>0$ (independent of $t_{0}$), there exists $T=T(\epsilon,c)>0$ such that, for all $t_{0}\geq0$, it holds
\begin{equation}
\mathbb{P}\big(|x(t,t_{0},x_{0})|<k, \forall t\geq t_{0}+T\big)\geq1-\epsilon\quad \text{for all} \
x_{0}\in\mathbb{R}^d \ \text{such that}\ |x_{0}|<c.
\end{equation}
\end{defn}

 Let us now denote by $C^{2,1}(\mathbb{R}^{d}\times\mathbb{R}_{+},\mathbb{R}_{+})$ the family of all nonnegative functions $V(x,t)$ defined on $\mathbb{R}^{d}\times\mathbb{R}_{+}$ which are twice continuously differentiable in $x$ and once in $t$.
Thanks to the It\^{o} formula it follows\\
$$\displaystyle{dV(x(t),t)=LV(x(t),t)dt+V_{x}(x(t),t)g(x(t),t)dW_{t}},$$
where
$$LV(x,t)=V_{t}(x,t)+V_{x}(x,t)f(x,t)+\frac{1}{2}trace[g(x,t)^{T}V_{xx}(x,t)g(x,t)],$$
$\displaystyle{V_{t}(x,t)=\frac{\partial V}{\partial t}(x,t)}$ ; $\displaystyle{V_{x}(x,t)=(\frac{\partial V}{\partial x_{1}}(x,t),...,\frac{\partial V}{\partial x_{d}}(x,t))}$;
$\displaystyle{V_{xx}(x,t)=\Big(\frac{\partial^{2} V}{\partial x_{i}\partial x_{j}}(x,t)\Big)_{d\times d}}.$

Let us also recall that by $\mathcal{K}$ it is usual to denote the family of all continuous nondecreasing functions $\mu$ : $\mathbb{R_{+}}\rightarrow\mathbb{R}_{+}$ such that $\mu(0)=0$, and by $\mathcal{K_{\infty}}$ the class of all functions $\mu\in\mathcal{K}$ such that\\ $\mu(r)\rightarrow+\infty$, $r\rightarrow+\infty$.

Now we can prove our first result.
\begin{thm}
Assume that there exist $V\in C^{2,1}(\mathbb{R}^{d}\times\mathbb{R}_{+},\mathbb{R}_{+})$, $\mu_{1}, \mu_{2}\in\mathcal{K}_{\infty}$, with\\ $\mu_{1}$ : $\mathbb{R}^{*}_{+}\rightarrow\mathbb{R}^{*}_{+}$ and $\mu_{3}\in\mathcal{K}$, $M>0$,  such that for all $t\geq0$, and all $x\in\mathbb{R}^{d}$,
$$\mu_{1}(|x|)\leq V(x,t)\leq\mu_{2}(|x|),$$

$$LV(x,t)\leq\varrho(t)-\mu_{3}(|x|),$$
where $\varrho(t)$ is a continuous nonnegative function with $$\lim_{t\rightarrow+\infty}\varrho(t)=0,$$ and $$\int_{0}^{+\infty}\varrho(t)dt\leq M<+\infty.$$
Then, system (\ref{la}) is globally uniformly practically stable in probability.
\end{thm}
\noindent {\bf{Proof.}}
Let $\epsilon\in]0,1[$, and take a small positive real number $r$,  and $k>r$. Assume that there exists $\delta=\delta(\epsilon,k)\in]0,k[$ such that, for $x_{0}\in S_{\delta}$, $\displaystyle{S_{\delta}=\{x\in\mathbb{R}^{d}/ |x|<\delta\}}$ and $\displaystyle{\mu_{1}(k)>\frac{M}{\epsilon}}$, we deduce
$$\frac{V(x_{0},0)}{\epsilon}\leq\mu_{1}(k)-\frac{M}{\epsilon}.$$
Let us write $x(t,t_{0},x_{0})=x(t)$ for simplicity. Define $\displaystyle{\tau=\inf\{t\geq0; |x(t)|\geq k\}}$ (Throughout this paper we set $\inf\emptyset=\infty$). By the It\^{o} formula, for any $t\geq0$, we have\\
$$V(x(\tau\wedge t),\tau\wedge t)=V(x_{0},0)+\int_{0}^{\tau\wedge t}LV(x(s),s)ds+\int_{0}^{\tau\wedge t}V_{x}(x(s),s)g(x(s),s)dW_{s}.$$
Then, $\displaystyle{0\leq\mathbb{E}\big(V(x(\tau\wedge t),\tau\wedge t)\big)=V(x_{0},0)+\mathbb{E}\big(\int_{0}^{\tau\wedge t}LV(x(s),s)ds\big)}$ and
\begin{align*}
\mathbb{E}\big(V(x(\tau\wedge t),\tau\wedge t)\big)&\leq V(x_{0},0)+\mathbb{E}\big(\int_{0}^{\tau\wedge t}\varrho(s)-\mu_{3}(|x(s)|)ds\big)\\
&\leq V(x_{0},0)+\mathbb{E}\big(\int_{0}^{\tau\wedge t}\varrho(s)ds\big)\\
&\leq V(x_{0},0)+\mathbb{E}\big(\int_{0}^{+\infty}\varrho(s)ds\big)\\
&\leq V(x_{0},0)+M.
\end{align*}
Notice that $\displaystyle{|x(\tau\wedge t)|=|x(\tau)|=k}$ if $\tau\leq t$.\\
Therefore, we can obtain
\begin{align*}
\mathbb{E}\big(V(x(\tau\wedge t),\tau\wedge t)\big)&\geq\mathbb{E}\big(1_{\{\tau\leq t\}}V(x(\tau),\tau)\big)\\
&\geq\mathbb{E}\big(1_{\{\tau\leq t\}}\mu_{1}(|x(\tau)|)\big)\\
&\geq\mathbb{E}\big(1_{\{\tau\leq t\}}\mu_{1}(k)\big)\\
&=\mu_{1}(k)\mathbb{E}\big(1_{\{\tau\leq t\}}\big)\\
&=\mu_{1}(k)\mathbb{P}\big(\tau\leq t\big).
\end{align*}
Then,
$$\displaystyle{\mu_{1}(k)\mathbb{P}\big(\tau\leq t\big)\leq V(x_{0},0)+M\leq\epsilon\mu_{1}(k)}.$$
Thus, $$\displaystyle{\mathbb{P}\big(\tau\leq t\big)\leq\epsilon}.$$
Letting $t\rightarrow+\infty$,  we arrive at
 $$\displaystyle{\mathbb{P}\big(\tau<+\infty\big)\leq\epsilon}.$$
 Thus, for all $\epsilon\in]0,1[$, $r\geq0$ and $k>r$, there exists $\delta=\delta(\epsilon,k)\in]0,k[$, such that
\begin{equation}
\mathbb{P}\big(|x(t)|<k, \forall t\geq t_{0}\geq0\big)\geq1-\epsilon,\quad \text{for all}\ |x_{0}|<\delta.
\end{equation}
In particular, for a small $r\in]0,\delta[$ and $|x_{0}|<r$, we can take $k=k(r)>0$, such that for all $t\geq t_{0}\geq0$, we have
$\displaystyle{|x(t)|<k}$,  a.s., for all $|x_{0}|<r$, and thus, $$\displaystyle{\sup_{t\geq t_{0}}|x(t)|<k}\quad \text{for all}\ |x_{0}|<r\quad a.s.$$ Therefore, $B_{r}$ is globally uniformly stable in probability.\\

Let us now prove the uniform attractivity of $B_{r}$. To this end, it is sufficient to prove that, for each $\epsilon\in]0,1[$, $h>r$ and $c>0,$  there exists $T=T(\epsilon,c)>0$ (independent of $t_{0}$), such that for all $ t\geq t_{0}+T$, we have
\begin{equation}
\mathbb{P}\big(|x(t)|< h\big)\geq1-\frac{\epsilon}{2}, \ \text{for all} \ |x_{0}|<c.
\end{equation}
Pick $\epsilon\in]0,1[$, $h>r$ large enough and $c>0$ (independent of $t_{0}$), such that

$$\displaystyle{|x_{0}|<c}, \ \displaystyle{\mu_{1}(h)>\frac{2M}{\epsilon}}\ \ \text{and}\ \ \displaystyle{\mu_{1}(h)-\frac{2M}{\epsilon}\geq\frac{2V(x_{0},0)}{\epsilon}}.$$
Suppose now that there exists $T=T(\epsilon,c)>0$ such that, for all $t_{0}\geq0$,
$$\displaystyle{\tau_{h}=\inf\{t\geq t_{0}+T; |x(t)|\geq h\}}.$$
Noting that, $\displaystyle{|x(\tau_{h}\wedge t)|=|x(\tau_{h})|=h}$ when $\tau_{h}\leq t$, and
applying the It\^{o} formula, we obtain, for all $t\geq t_{0}+T$,
\begin{align*}
\mathbb{E}\big(V(x(\tau_{h}\wedge t),\tau_{h}\wedge t)\big)&\leq V(x_{0},0)+\mathbb{E}\big(\int_{0}^{\tau_{h}\wedge t}(\varrho(s)-\mu_{3}(|x(s)|)ds\big),\\
&\leq V(x_{0},0)+M,\\
&\leq\frac{\epsilon}{2}\mu_{1}(h)-M+M,\\
&\leq\frac{\epsilon}{2}\mu_{1}(h).
\end{align*}
We also have
\begin{align*}
\mathbb{E}\big(V(x(\tau_{h}\wedge t),\tau_{h}\wedge t)\big)&\geq\mathbb{E}\big(1_{\{\tau_{h}\leq t\}}V(x(\tau_{h}),\tau_{h})\big)\\
&\geq\mathbb{E}\big(1_{\{\tau_{h}\leq t\}}\mu_{1}(|x(\tau_{h})|)\big)\\
&\geq\mathbb{E}\big(1_{\{\tau_{h}\leq t\}}\mu_{1}(h)\big)\\
&=\mu_{1}(h)\mathbb{E}\big(1_{\{\tau_{h}\leq t\}}\big)\\
&=\mu_{1}(h)\mathbb{P}\big(\tau_{h}\leq t\big).
\end{align*}
Therefore,
$$\displaystyle{\mu_{1}(h)\mathbb{P}\big(\tau_{h}\leq t\big)\leq\frac{\epsilon}{2}\mu_{1}(h)},$$
which implies that $$\displaystyle{\mathbb{P}\big(\tau_{h}\leq t\big)\leq\frac{\epsilon}{2}}.$$
Then, for any $\epsilon\in]0,1[$, $h>r$ and $c>0$ (independent of $t_{0}$), for a fixed $T=T(\epsilon,c)>0$, we have, for all $t\geq t_{0}+T$,
$\displaystyle{\mathbb{P}\big(\tau_{h}\leq t\big)\leq\frac{\epsilon}{2}}$, for all $|x_{0}|<c$. Letting $t\rightarrow+\infty$, we arrive at
 $$\displaystyle{\mathbb{P}\big(\tau_{h}<+\infty\big)\leq\frac{\epsilon}{2}}.$$ Thus,
$$\mathbb{P}\big(|x(t)|<h; \forall t\geq t_{0}+T \big)\geq1-\frac{\epsilon}{2},\quad \text{for all}\ |x_{0}|<c.$$
Then, $B_{r}$ is  globally uniformly attractive in probability, as desired.
\hfill $\Box$
\begin{rem}
Note that if $r\rightarrow0$ or $\varrho(t)\rightarrow0$ when $t\rightarrow+\infty$,   we then recover the classical stability results when $0$ is an equilibrium point.
\end{rem}
\section{The global practical uniform exponential stability of stochastic differential equation}
Let us start by recalling a technical lemma which will be useful in our analysis.
\begin{lem}\label{ar}
For all $x_{0}\in\mathbb{R}^{d}$ such that $x_0\neq0$ it holds
\begin{equation}
\displaystyle{\mathbb{P}\big(x(t,t_{0},x_{0})\neq0,\forall t\geq0\big)=1}.
\end{equation}
\end{lem}
\noindent {\bf{Proof.}}
Arguing by contradiction, if the lemma were false, there would exist some $x_{0}\neq0$ such that $\displaystyle{\mathbb{P}\{\tau<\infty\}>0}$, where
$$\tau=\inf\{t\geq0; x(t)=0\},$$
where we write $x(t,t_{0},x_{0})=x(t)$ for simplicity. Thus, we can find a pair of constants $T>0$ and $\theta>1$, sufficiently large, such that $\mathbb{P}(B)>0$, where
$$B=\{\tau\leq T\quad \text{and}\quad |x(t)|\leq\theta-1,\quad \text{for all}\quad 0\leq t\leq\tau\}.$$
But, by the standing hypotheses, there exist positive constants $K_{\theta}$, $r$, and a small positive constant $c$ such that
$$|f(x,t)|\vee|g(x,t)|\leq K_{\theta}|x|+r \quad \text{for all}\ 0<c\leq|x|\leq\theta,\quad 0\leq t\leq T.$$
Let $\displaystyle{V(x,t)=|x|^{-1}}$. Then, for $0<c\leq|x|\leq\theta$ and $0\leq t\leq T$,
\begin{align*}
LV(x,t)&=-|x|^{-3}x^{T}f(x,t),\\
&\qquad +\frac{1}{2}\Big(-|x|^{-3}|g(x,t)|^{2}+3|x|^{-5}|x^{T}g(x,t)|^{2}\Big)\\
&\leq|x|^{-2}|f(x,t)|+|x|^{-3}|g(x,t)|^{2}\\
&\leq|x|^{-2}\big(K_{\theta}|x|+r\big)+|x|^{-3}\big(K_{\theta}|x|+r\big)^{2}\\
&\leq K_{\theta}(1+K_{\theta})V(x,t)+r|x|^{-1}\big(|x|^{-1}+2K_{\theta}|x|^{-1}+r|x|^{-2}\big)\\
&\leq K_{\theta}(1+K_{\theta})V(x,t)+rV(x,t)\big(\frac{1}{c}+\frac{2K_{\theta}}{c}+\frac{r}{c^{2}}\big)\\
&\leq \big[K_{\theta}(1+K_{\theta})+\beta\big]V(x,t),
\end{align*}
where $\displaystyle{\beta=r\big(\frac{1}{c}+\frac{2K_{\theta}}{c}+\frac{r}{c^{2}}\big)}$.

Now, for any $\epsilon\in]0,|x_{0}|[$, define the stopping time
$$\tau_{\epsilon}=\inf\{t\geq0; |x(t)|\in]-\infty,\epsilon]\cup[\theta,+\infty[\}.$$
Again by the It\^{o} formula,
\begin{align*}
&  \mathbb{E}\Big[\exp\Big\{-\big[K_{\theta}(1+K_{\theta})+\beta\big](\tau_{\epsilon}\wedge T)\Big\}V\big(x(\tau_{\epsilon}\wedge T),\tau_{\epsilon}\wedge T\big)\Big]-V(x_{0},0),\\
&=\mathbb{E}\Big[\int_{0}^{\tau_{\epsilon}\wedge T}\exp\Big\{-\big[K_{\theta}(1+K_{\theta})+\beta\big]s\Big\}\Big[-\big(K_{\theta}(1+K_{\theta})+\beta\big)V(x(s),s)+LV(x(s),s)\Big]ds,\\
&\leq0.
\end{align*}
Note that for $\omega\in B$, it follows that $\tau_{\epsilon}\leq T$ and $|x(\tau_{\epsilon})|=\epsilon$. The above inequality therefore implies that
$$\mathbb{E}\Big[\exp\Big\{-\big[K_{\theta}(1+K_{\theta})+\beta\big]T\Big\}\epsilon^{-1}1_{B}\Big]\leq|x_{0}|^{-1}.$$
Hence $\displaystyle{\mathbb{P}(B)\leq\epsilon|x_{0}|^{-1}\exp\Big\{\big[K_{\theta}(1+K_{\theta})+\beta\big]T\Big\}}$. Letting $\epsilon\longrightarrow0$ yields that $\displaystyle{\mathbb{P}(B)=0}$, but this contradicts the definition of $B$. The proof is complete.
\hfill $\Box$\\

The next lemma, whose proof can be found in \cite{13}, will be crucial in the proof of our main result.
\begin{lem}\label{la2}
Let $\displaystyle{g=(g_{1},....,g_{m})\in L^{2}(\mathbb{R}_{+},\mathbb{R}^{m})}$, and let  $T$, $\alpha$, $\beta$ be any positive numbers. Then
$$\mathbb{P}\Big(\sup_{0\leq t\leq T}\big[\int_{0}^{t}g(s)dW_{s}-\frac{\alpha}{2}\int_{0}^{t}|g(s)|^{2}ds\big]>\beta\Big)\leq\exp(-\alpha\beta).$$
\end{lem}
\begin{defn}
The ball $B_r$ is said to be almost surely globally practically uniformly exponentially stable if, for any $x_0$ such that $0<|x(t,t_0,x_0)|-r$, for all  $t\geq t_0\geq0$, it holds that
\begin{equation}\label{re}
\lim_{t\to\infty}\sup \frac{1}{t}\ln(|x(t,t_0,x_0)|-r)<0, \ \textrm{a.s.}
\end{equation}

System (\ref{la}) is said to be almost surely globally practically uniformly exponentially stable  if there exists $r>0$ such that $B_{r}$ is almost surely globally practically uniformly
 exponentially stable.
\end{defn}
Now, we can establish and prove our main result in this paper.
\begin{thm}\label{la4}
Assume that there exist a function $V\in C^{2,1}(\mathbb{R}^{d}\times\mathbb{R}_{+},\mathbb{R^{*}}_{+})$ and constants $p\in\mathbb{N}^{*}$, $c_{1}\geq1$, $\varrho\geq c_{1}$, $\gamma\geq0$ and $c_{2}\in\mathbb{R}$, $c_{3}\geq0$ such that for all $ t\geq t_{0}\geq0$, and $x\in\mathbb{R}^d$, the following conditions hold:

$$c_{1}|x|^{p}\leq V(x,t),$$
$$LV(x,t)\leq c_{2}V(x,t)+\varrho,$$
$$|V_{x}(x,t)g(x,t)|^{2}\geq c_{3}V^{2}(x,t)+\gamma.$$
Then
$$\displaystyle{\lim_{t\rightarrow+\infty}\sup\frac{1}{t}\ln\big(|x(t,t_{0},x_{0})|-(\frac{\varrho}{c_{1}})^{\frac{1}{p}}\big)\leq -\frac{\big[c_{3}-2(c_{2}+1)\big]}{2}},\quad a.s.,\ \text{for all}\ x_{0}\in\mathbb{R}^{d}$$
 In particular, if $\displaystyle{c_{3}>2(c_{2}+1)}$, then the solution of (\ref{la}) is almost surely  globally practically uniformly exponentially stable with $\displaystyle{r=(\frac{\varrho}{c_{1}})^{\frac{1}{p}}}$.
\end{thm}
\noindent {\bf{Proof.}}
Fix $x_{0}\neq0$ in $\mathbb{R}^{d}$ and write $x(t,t_{0},x_{0})=x(t)$. By Lemma \ref{ar}, $x(t)\neq 0$, for all $t\geq0$ almost surely.

Observe that, for  any $x\in\mathbb{R}^d$, we have
\begin{align*}
c_{1}|x|^{p}-\varrho&=c_{1}\big(|x|^{p}-\frac{\varrho}{c_{1}}\big),\\
&=c_{1}\big(|x|^{p}-\big((\frac{\varrho}{c_{1}})^{\frac{1}{p}}\big)^{p}\big),\\
&=c_{1}\big(|x|-(\frac{\varrho}{c_{1}})^{\frac{1}{p}}\big)\big(|x|^{p-1}+|x|^{p-2}(\frac{\varrho}{c_{1}})^{\frac{1}{p}}+...+(\frac{\varrho}{c_{1}})^{\frac{p-1}{p}}\big),\\
&\geq c_{1}\big(|x|-(\frac{\varrho}{c_{1}})^{\frac{1}{p}}\big).
\end{align*}
Since $c_{1}\geq1$, we obtain
$$V(x,t)\geq  c_{1}|x|^{p}\geq c_{1}|x|^{p}-\varrho\geq\big(|x|-(\frac{\varrho}{c_{1}})^{\frac{1}{p}}\big).$$

Therefore, $\big(|x|-(\frac{\varrho}{c_{1}})^{\frac{1}{p}}\big)\leq V(x,t)$ and $\ln\big(|x|-(\frac{\varrho}{c_{1}})^{\frac{1}{p}}\big)\leq\ln\big(V(x,t)\big)$,  for all $t\geq0, x\in\mathbb{R}^d$.
Applying the It\^{o} formula once more and taking into account the assumptions, we obtain that, for all $t\geq0$
\begin{align*}
d\big(\ln(V(x(t),t))\big)&=\frac{LV(x(t),t)}{V(x(t),t)}dt+\frac{V_{x}(x(t),t)g(x(t),t)}{V(x(t),t)}dW_{t},\\
&-\frac{1}{2}\frac{|V_{x}(x(t),t)g(x(t),t)|^{2}}{V^{2}(x(t),t)}dt,
\end{align*}
\begin{align*}
\int_{0}^{t}d\big(\ln(V(x(s),s))\big)ds&=\int_{0}^{t}\frac{LV(x(s),s)}{V(x(s),s)}ds+\int_{0}^{t}\frac{V_{x}(x(s),s)g(x(s),s)}{V(x(s),s)}dW_{s},\\
&-\frac{1}{2}\int_{0}^{t}\frac{|V_{x}(x(s),s)g(x(s),s)|^{2}}{V^{2}(x(s),s)}ds,
\end{align*}
\begin{eqnarray*}
\ln(V(x(t),t))&=&\ln(V(x(0),0))+\int_{0}^{t}\frac{LV(x(s),s)}{V(x(s),s)}ds+M(t),\\
&&-\frac{1}{2}\int_{0}^{t}\frac{|V_{x}(x(s),s)g(x(s),s)|^{2}}{V^{2}(x(s),s)}ds,\\
&\leq&\ln(V(x(0),0))+\int_{0}^{t}\frac{c_{2}V(x(s),s)+\varrho}{V(x(s),s)}ds+M(t),\\
&&-\frac{1}{2}\int_{0}^{t}\frac{|V_{x}(x(s),s)g(x(s),s)|^{2}}{V^{2}(x(s),s)}ds,\\
&\leq&\ln(V(x(0),0))+c_{2}t+\int_{0}^{t}\frac{\varrho}{V(x(s),s)}ds+M(t),\\
&&-\frac{1}{2}\int_{0}^{t}\frac{|V_{x}(x(s),s)g(x(s),s)|^{2}}{V^{2}(x(s),s)}ds,\\
&\leq&\ln(V(x(0),0))+c_{2}t+\int_{0}^{t}\frac{\varrho}{c_{1}|x(s)|^{p}}ds+M(t),\\
&&-\frac{1}{2}\int_{0}^{t}\frac{|V_{x}(x(s),s)g(x(s),s)|^{2}}{V^{2}(x(s),s)}ds,\\
&\leq&\ln(V(x(0),0))+c_{2}t+t+M(t)-\frac{1}{2}\int_{0}^{t}\frac{|V_{x}(x(s),s)g(x(s),s)|^{2}}{V^{2}(x(s),s)}ds,
\end{eqnarray*}
\begin{equation}\label{la3}
\ln(V(x(t),t))\leq\ln(V(x(0),0))+(c_{2}+1)t+M(t)-\frac{1}{2}\int_{0}^{t}\frac{|V_{x}(x(s),s)g(x(s),s)|^{2}}{V^{2}(x(s),s)}ds,
\end{equation}
where $\displaystyle{M(t)=\int_{0}^{t}\frac{V_{x}(x(s),s)g(x(s),s)}{V(x(s),s)}dW_{s}}$ is a continuous martingale with initial value\\ $M(0)=0$. Assign $\epsilon\in]0,1[$ arbitrarily and let  $n=1,2,...$ By Lemma \ref{la2},
$$\mathbb{P}\Big\{\sup_{0\leq t\leq n}\big[M(t)-\frac{\epsilon}{2}\int_{0}^{t}\frac{|V_{x}(x(s),s)g(x(s),s)|^{2}}{V^{2}(x(s),s)}ds\big]>\frac{2}{\epsilon}\ln(n)\Big\}\leq\frac{1}{n^{2}}.$$
Applying the Borel-Cantelli lemma we see that, for almost all $\omega\in\Omega$, there exists an integer\\ $n_{0}=n_{0}(\omega)$, such that if $n\geq n_{0}$,
$$M(t)\leq\frac{2}{\epsilon}\ln(n)+\frac{\epsilon}{2}\int_{0}^{t}\frac{|V_{x}(x(s),s)g(x(s),s)|^{2}}{V^{2}(x(s),s)}ds,\ \text{for all}\  0\leq t\leq n.$$
Then, since $\displaystyle{\frac{\epsilon-1}{2}\int_{0}^{t}\frac{\gamma}{V^{2}(x(s),s)}ds\leq0}$,  inequality (\ref{la3}) becomes, for all $0\leq t\leq n$, $n\geq n_{0}$ almost surely, as
\begin{align*}
\ln(V(x(t),t))&\leq\ln(V(x(0),0))+(c_{2}+1)t+\frac{2}{\epsilon}\ln(n)+\frac{\epsilon-1}{2}\int_{0}^{t}\frac{|V_{x}(x(s),s)g(x(s),s)|^{2}}{V^{2}(x(s),s)}ds,\\
&\leq\ln(V(x(0),0))+(c_{2}+1)t+\frac{2}{\epsilon}\ln(n)+\frac{\epsilon-1}{2}c_{3}t+\frac{\epsilon-1}{2}\int_{0}^{t}\frac{\gamma}{V^{2}(x(s),s)}ds,\\
&\leq\ln(V(x(0),0))+\frac{2}{\epsilon}\ln(n)-\frac{1}{2}\big[(1-\epsilon)c_{3}-2(c_{2}+1)\big]t.
\end{align*}
Consequently, for almost all $\omega\in\Omega$, if $n-1\leq t\leq n$ and $n\geq n_{0}$, we deduce
$$\frac{1}{t}\ln(V(x(t),t))\leq-\frac{1}{2}\big[(1-\epsilon)c_{3}-2(c_{2}+1)\big]+\frac{\ln(V(x(0),0))+\frac{2}{\epsilon}\ln(n)}{n-1}.$$
This implies that,
$$\lim_{t\rightarrow+\infty}\sup\frac{1}{t}\ln(V(x(t),t))\leq-\frac{1}{2}\big[(1-\epsilon)c_{3}-2(c_{2}+1)\big]\quad \text{a.s.}$$
Hence,
$$\lim_{t\rightarrow+\infty}\sup\frac{1}{t}\ln(|x(t)|-(\frac{\varrho}{c_{1}})^{\frac{1}{p}})\leq-\frac{1}{2}\big[(1-\epsilon)c_{3}-2(c_{2}+1)\big]\quad \text{a.s.}$$
and the required assertion follows since $\epsilon>0$ is arbitrary.
\hfill $\Box$
\begin{rem}
Note that, as the origin $x=0$ may not be an equilibrium point of the system (\ref{la}), then we can no longer study the stability of the origin as an equilibrium point nor should we expect the solution of the system to approach the origin almost surely as $t\rightarrow+\infty$. The inequality (\ref{re}) implies that $x(t)$ will be ultimately bounded by a small bound $r>0$, that is, $|x(t)|$ will be small for sufficiently large $t$. This can be viewed as a robustness property of convergence almost surely to the origin provided that $f$ and $g$ satisfies $f(0,t)=0$ and $g(0,t)=0$, $\forall t\geq0$. In this case the origin becomes an equilibrium point.
\end{rem}
\begin{exmp}
Although our theory may be applied to a more general situation, we will exhibit now a very simple situation in order to illustrate how our main result works. Indeed, let us consider the Langevin equation (see \cite{5}, \cite{14}):
\begin{equation}\label{la20}
dx(t)=f(x(t))dt+g(x(t))dW_{t}=\alpha x(t)dt+\beta dW_{t},
\end{equation}
where $W_{t}$ is a one dimensional standard Brownian motion and $f(x,t)=\alpha x$, $g(x,t)=\beta$, $V(x,t)=x^{2}$, $\beta\in\mathbb{R}$, $|x|>\sqrt{\beta^{2}+1}$ and $\alpha<-\frac{1}{2}$. Then, it is easy to check that
$$LV(x,t)=V_{t}(x,t)+V_{x}(x,t)f(x,t)+\frac{1}{2}g(x,t)^{2}V_{xx}(x,t)=2\alpha x^{2}+\beta^{2}.$$
$$|V_{x}(x,t)g(x,t)|^{2}=|2\beta x|^{2}=4\beta^{2} x^{2}.$$
Therefore,
$$|x|^{2}\leq V(x,t),$$
$$LV(x,t)\leq2\alpha V(x,t)+\beta^{2}+1,$$
$$|V_{x}(x,t)g(x,t)|^{2}=4\beta^{2} x^{2}\geq0.$$
Thus, the constants in Theorem \ref{la4} become  $c_{1}=1$, $c_{2}=2\alpha$, $c_{3}=0$, $p=2$, $\varrho=\beta^{2}+1$, $\gamma=0$.\\

Consequently, $$c_{3}=0>2(c_{2}+1)=4\alpha+2,$$
and the assumptions of Theorem \ref{la4} are fulfilled. Then, (\ref{la20}) is almost sure globally practically uniformly exponentially stable with $r=\sqrt{\beta^{2}+1}$.
\end{exmp}
\begin{rem}
\begin{enumerate}
\item As we mentioned in the Introduction,  when the origin is not a solution of the stochastic equation, the practical stability means that almost all the paths of each solution starting in deterministic initial values, remain inside a neighborhood of certain  ball centered at the origin. This definition is therefore weaker than the usual one of almost sure stability, in which the initial values are allowed to be non-deterministic.
\item On the other hand, and in some sense, one could interpret the global practical stability established in Theorem \ref{la4} as a result ensuring the existence of an absorbing set for the solutions starting in deterministic initial values. Therefore, it would be very interesting to analyze whether this kind of results are related to the theory of random attractors (see e.g. \cite{Schmalfuss92}, \cite{Kloeden}, \cite{crauel-flandoli}). However, as the forward convergence stated in Theorem \ref{la4} is not necessarily uniform it is not clear that we might ensure the existence of a forward random attractor. 

\item Nevertheless, in the previous example it is known the existence of a pullback random attractor given by the so-called Ornstein-Uhlenbeck process, which is a special solution of the equation \eqref{la20}, corresponding to a non-deterministic initial value, and which attracts any other solution of the equation in the pullback  sense (see \cite{Kloeden} for more details).   

\end{enumerate}
\end{rem}

\textbf{Acknowledgements.} We would like to thank the referee for the interesting comments and suggestions which allowed us to improve the presentation of this paper.

\providecommand{\bysame}{\leavevmode\hbox
to3em{\hrulefill}\thinspace}
\providecommand{\MR}{\relax\ifhmode\unskip\space\fi MR }
\providecommand{\MRhref}[2]{%
  \href{http://www.ams.org/mathscinet-getitem?mr=#1}{#2}
} \providecommand{\href}[2]{#2}

\end{document}